\documentclass[12pt]{article}

\usepackage{latexsym,amsmath,amssymb}
\usepackage{graphicx}

\title{\bfseries%
  Two classes of quadratic vector fields for which the Kahan discretization is integrable
}
\author{%
  Elena Celledoni$^1$, Robert I. McLachlan $^2$, David I. McLaren$^3$, \\
  Brynjulf Owren$^1$, G.R.W. Quispel$^3$,\\[12pt]
  \normalsize
  1~Norwegian University of Science and Technology, 7491, Trondheim, Norway\\
  \normalsize
  \normalsize
  2~Massey University, Palmerston North, New Zealand\\
\normalsize
  3~La Trobe University, Victoria 3086, Australia\\
  \normalsize
}
\date{}

\begin{document}

\maketitle
\pagestyle{empty}
\thispagestyle{empty}

\section{INTRODUCTION}
\label{sec:INTRODUCTION}

Kahan's method for discretizing quadratic differential equations was introduced in \cite{Kahan}. It was rediscovered in the context of integrable systems by Hirota and Kimura \cite{HirKim}. Suris and collaborators extended the applications to integrable systems significantly in a series of papers \cite{PetSur1}, \cite{PetSur2}, \cite{PetSur3}, \cite{PetPfaSur}, \cite{HonePet}. Applications to non-integrable Hamiltonian systems and the use of polarisation to discretise arbitrary degree Hamiltonian systems were studied in \cite{celledoni1}, \cite{celledoni2} and \cite{celledoni3}.

The present paper contains an extract of the talk one of the authors (GRWQ) gave on 6th July 2016 at the 12th International Conference on Symmetries and Integrability of Difference Equations (SIDE12) in Sainte Adele, Quebec, Canada. 
We present two classes of $2$-dimensional ODE systems of quadratic vector fields where the Kahan discretization is integrable. Both classes of systems can be cast  in the form
 \begin{equation}
 \label{eq:1}
 \frac{dX}{dt} = \varphi(X) K \nabla H(X),
 \end{equation}
 where 
 $$X^t=(x,y),\qquad K= \left( \begin{array}{cc}
0 & 1  \\
-1 & 0  \end{array} \right),$$
and $\varphi(X)$ is a scalar function of the components of $X$.
In the first class the Hamiltonian function is quartic and in the second class it is sextic. These systems can be seen as generalisations of the examples of the reduced Nahm equations presented in \cite{PetPfaSur}.
Some of the results in this paper were  found independently by Petrera and Zander \cite{PetZan}.

\section{Quartic Hamiltonians in $2$D} 

Consider the $2$-dimensional ODE system \eqref{eq:1} where $\varphi(X)=\frac{1}{ax+by}$,
and the homogeneous Hamiltonian has the form $H=(ax+by)^2(cx^2+2dxy+ey^2)$. Then the Kahan map for this system preserves the modified Hamiltonian:
{\footnotesize
\begin{align*}
\widetilde{H}(X)=\frac{H}{(1+h^2D(ax+by)^2
+h^2E(cx^2+2dxy+ey^2))(1+h^{2}9D(ax+by)^2
+h^2E(cx^2+2dxy+ey^2)},
\end{align*}
}
\small{
and the measure:
\begin{align*}
m(x,y)= \frac{dxdy}{(ax+by)(cx^2+2dxy+ey^2)}.
\end{align*}
Here, $D:=ce-d^2$ and $E:=2abd-a^2e-b^2c$. It follows that the Kahan map is integrable. It also seems to preserve the genus of the level sets of $H$ (which generically equals 1).

\section{Sextic Hamiltonians in $2$D} 

Consider the $2$-dimensional ODE system \eqref{eq:1} with $\varphi(X)=\frac{1}{(cx+dy)(ex+fy)^2}$ and with the {\em homogeneous} sextic 
Hamiltonian $H=(ax+by)(cx+dy)^2(ex+fy)^3$. Then the Kahan map for this system preserves the modified Hamiltonian:
\small{
\begin{align*}
\widetilde{H}(X)=\frac{H}{(1+a_{5}l_{2}^2)(1+a_{3}l_{1}^2+a_{4}l_{3}^2+a_{7}l_{1}l_{3})(1+a_{5}l_{2}^2+a_{6}l_{3}^2))}
\end{align*}
} where
\begin{align*}
l_{1}:=ax+by, \quad  l_{2}:=cx+dy, \quad l_{3}:=ex+fy,
\end{align*}
and $d_{1,2}:=h(ad-bc)$, $d_{2,3}:=h(cf-ed)$, $d_{3,1}:=h(eb-fa)$
and $a_{3}:=\frac{-9d_{2,3}^2}{4}$, $a_{4}:=\frac{-d_{1,2}^2}{4}$,
$a_{5}:=\frac{-9d_{3,1}^2}{4}$; $a_{6}:=-4d_{1,2}^2$, $a_{7}:=\frac{3d_{1,2}d_{2,3}}{2}$.
In this case the Kahan map also preserves  the modified measure
\begin{align*}
m(x,y)= \frac{dxdy}{(ax+by)(cx+dy)(ex+fy)}.
\end{align*}
Again, the Kahan map is integrable. It also seems to preserve the genus of the level sets of $H$ (which  generically equals 1).

%
\section*{Acknowledgment}

This work was supported by the Australian Research Council, by the Research Council of Norway, by the Marsden Fund of the Royal Society of New Zealand,
and by the European Union's Horizon 2020
research and innovation programme under the Marie Sklodowska-Curie grant
agreement No. 691070.


\end{document}